# Methods toward better Multigrid Solver Convergence

by

**John T. Wallis**

Applied Mathematics
Research Report

DEPARTMENT OF MATHEMATICAL SCIENCES

KAIST 한국과학기술원
Korea Advanced Institute of Science and Technology

# Methods toward better Multigrid Solver Convergence

John T. Wallis

Department of Mathematical Sciences
KAIST, Daejeon, 305-701, Korea

**Abstract**

*I present a motivation of several areas where the Multigrid techniques can be employed. I present typical areas where the multigrid solver might be employed. I give an introduction to smoothers and how one might choose a preconditionor as well as an introduction of the Multigrid technique used. Then I do a study of the Multigrid technique while adjusting the environment conditions of the solver and the problem such as anisotropies, the grid-levels used, the preconditionor smoothing steps, the coordinate system and the start vector. The problem solved here was a simple Poisson problem. The Multigrid program used an F-cycle in this paper. I include performance study sections displaying results of the solver behavior under the different conditions.*

**Introduction / Motivation**

Many iterative techniques employed to solve problems that arise from partial differential equations posses a smoothing property that reduces high frequency errors rapidly while low frequency errors are reduced slowly. Techniques to speed converges leads to methods that make a better start vector "guess" and also making use of multiple grids that aid in reducing low frequency errors. To improve multigrid iteration convergence, preconditionor techniques are employed with a parameter.

**Some Typical Problems:**
Problems are often formulated in different coordinate systems. Star simulations in astrophysics are often done using cylindrical or spherical coordinates. These are selected because they are close to the geometry. These geometries lead to complex coefficients that result in anisotropies that the solvers must handle.

Real life situations can have complex geometries. To accurately model their small scale geometries sufficiently special meshes must be used that can have high anisotropies. Again, solvers must be robust to handle these situations.

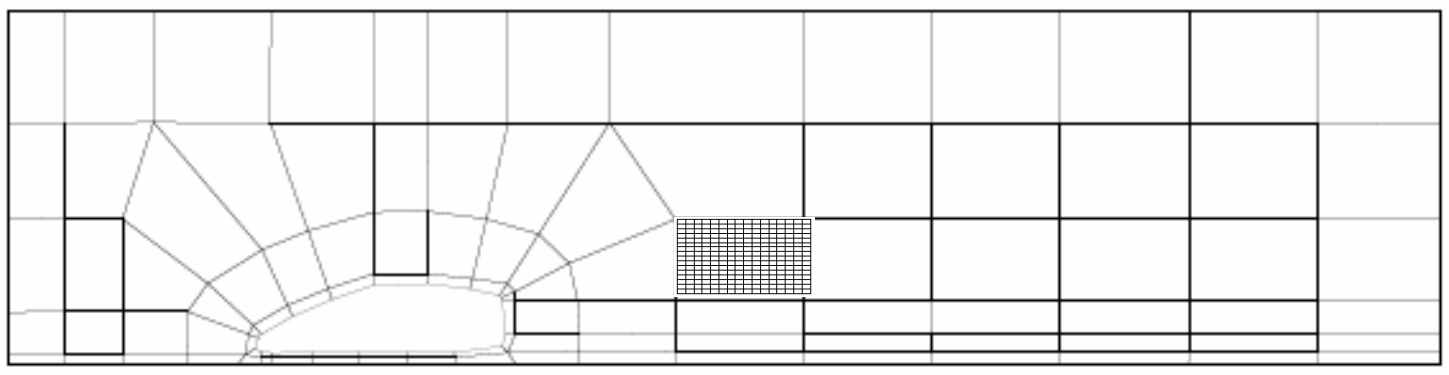
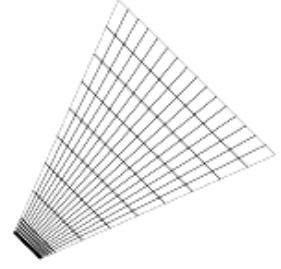
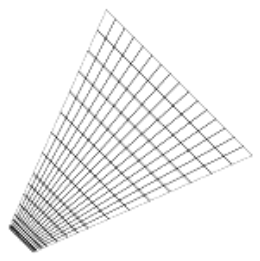

Problems can be large. Time and memory must be considered. Solvers must solve millions of unknowns on meshes of varying anisotropies with a high level of accuracy.

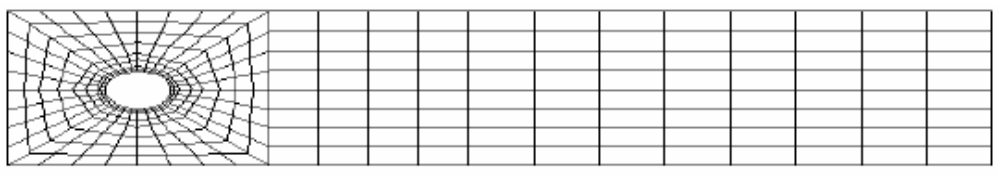
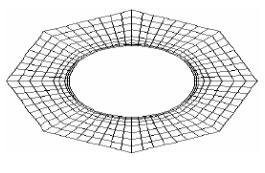
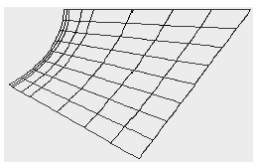

Matrix – vector operations take 60-90% of the computational cost. To reduce this sparse matrix techniques are employed. Also considerations are made to the PC-machine used. These include unrolling optimizations.

**In astronomy,** computations include the Euler/ momentum equation

$$\frac{\partial v}{\partial t} + v\frac{\partial v}{\partial x} + \frac{1}{\rho}\frac{dp}{dx} + \nabla\Phi = 0,$$

continuity equation $\frac{d\rho}{dt} + \frac{\partial}{\partial x}(\rho v) = 0$,as well as the Poisson equation $\Delta\Phi = -4\pi G\rho$ for the gravity potential $\Phi$. Such Laplace/ Poisson problems are the focus of my work.

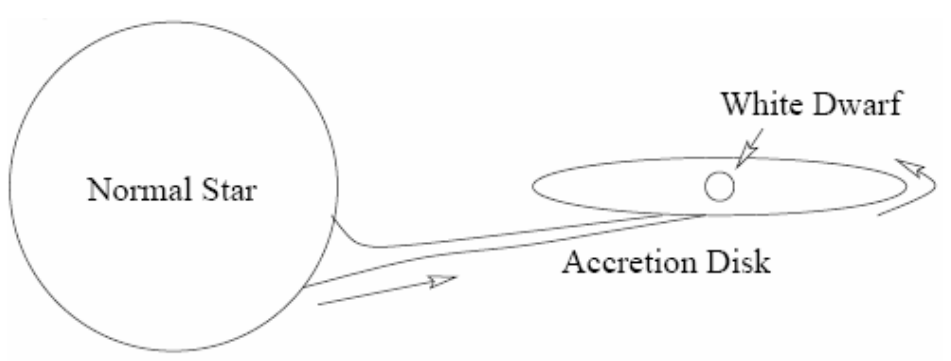

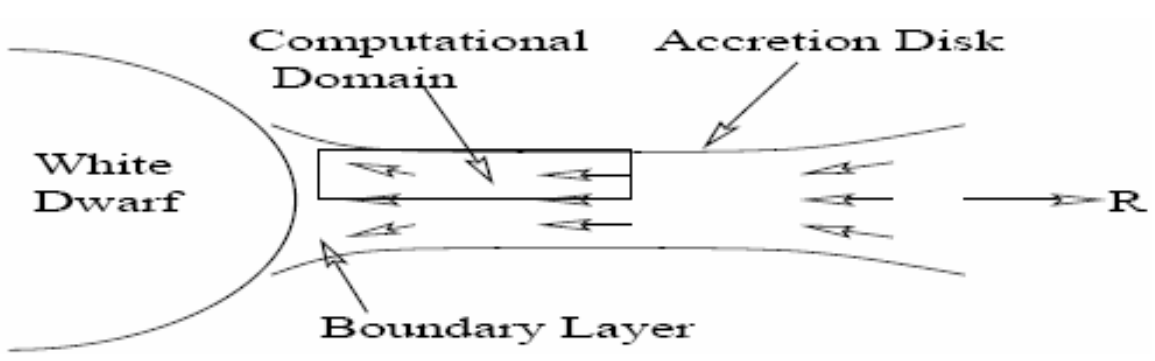

Below to the left is the grid of the problem and to the right its trans formation to the unit square

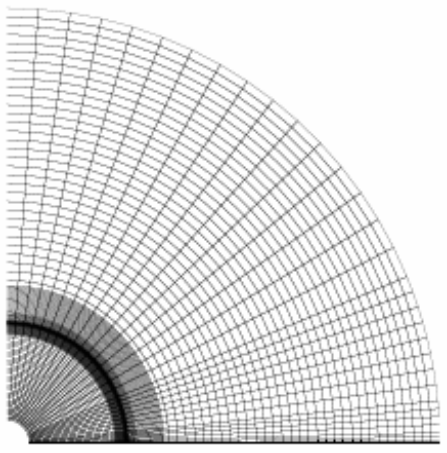 $\rightarrow$ 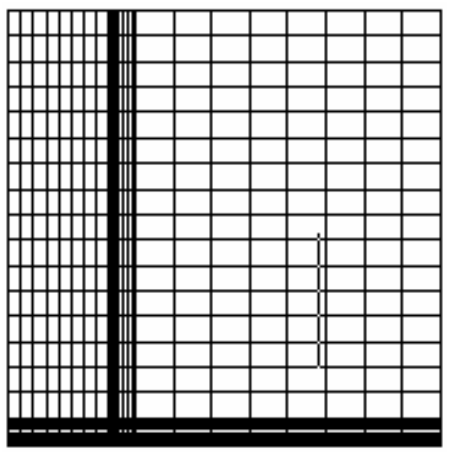

**A** popular problem is the Navier-Stokes (NS) formulation in computational fluid dynamics (CFD). It contains all the complexities and difficulties that make it ideal for study.
The incompressible NS problem is as follows:

$$\left.\begin{array}{l} u_t - \nu\Delta u + u\cdot\nabla u + \nabla p = f \\ \nabla\cdot u = 0 \end{array}\right\} in \quad \Omega\times(0,T].$$

with velocity u and pressure p. So $u_t$ represents the time derivative of u, $\nabla\cdot$ the divergence, $\Delta$ the Laplace operator and $\nabla$ the gradient.

A typical strategy to solve this would look as

$$\frac{u-u^n}{k} + \theta[-\nu\Delta u + u\cdot\nabla u] + \nabla p = g^{n+1}$$

$$with \quad r.h.s. \qquad g^{n+1} = \theta f^{n+1} + (1-\theta)f^n - (1-\theta)[-\nu\Delta u^n + u^n\cdot\nabla u^n]$$

After space and time discretization the corresponding problem can be written as

$$\begin{pmatrix} S & kB \\ B^T & 0 \end{pmatrix}\begin{pmatrix} u^L \\ p^L \end{pmatrix} = \begin{pmatrix} g \\ 0 \end{pmatrix},$$

where S is the coefficient matrix for u: Su = $[M+\theta kN(u)]u$ , $\theta$ = 1 (Backward Euler), $\theta$ = ½ (Crank Nicholson) and N(u)u = $-\nu\Delta u + u\cdot\nabla u$ = diffusive and convective part and B is the gradient and $B^T$ is the divergence.

A typical strategy is:

- Step 1: $u^L = S^{-1}(g - kBp^L)$ . Using the fact $B^T$ u = 0 Leads to $B^T u^L = B^T S^{-1}(g - kBp^L) = 0$
- Step 2: $f_p = \frac{1}{k}B^T u^L = \frac{1}{k}B^T S^{-1}(g - kBp^{L-1})$ =residual for $p^{L-1}$
- Step 3: Solve $-\Delta_h q = \frac{1}{k}B^T u^L$ Pressure-Poisson problem $-\Delta_h q = \frac{1}{k}B^T u^L$
- Step 4: Update $p^L = p^{L-1} + \alpha_R q + \alpha_D M_p^{-1} f_p$ , where $\alpha_R$, $\alpha_D$ are update parameters
- Step 5: Update $u^{L+1} = u^L - kBq$

**Smoothers**

The idea behind the preconditioned Richardson[1] iteration technique is

$$Ax = b \quad\Leftrightarrow\quad \omega Cx = \omega Cx - Ax + b \quad\Leftrightarrow\quad x = \quad x - \omega C^{-1}(Ax - b)$$

where C is an easy matrix to invert, good conditioned matrix that acts on the defect. For the sequence of solutions we get

$$x^{t+1} = x^t - \omega C^{-1}(Ax^t - b)$$
$$= (1 - C^{-1}A)x^t + \omega C^{-1}b =: Bx^t + c$$

So that the error looks like

$$e^{t+1} := x^{t+1} - x = Bx^t + x - (Bx + c) = Be^t$$

Of course, for convergence

$$\|B\| = \|I - \omega C^{-1}A\| < 1$$

So a clever choice of C will affect the behavior of the convergence.

The error vector is $e^t = x - x^t$ and the residual vector $r^t = Bx^t - b$ satisfies the relation $Be^{(n)} = r^t$. To accelerate convergence[2] the Multigrid technique then solves on a coarser level $\tilde{B}\tilde{e}^t = \tilde{r}^t$ and adds the prolonged result to the next finer grid level. The hope is we subtract the exact error: $x^{t+1} = x^t - e^t$

**The choice of C**

C is chosen between $\omega C^{-1} = A^{-1}$ and $C$ easy to invert. It is dependent on factors as computation cost and memory. Using D, L and U to denote the diagonal, lower and upper parts of A respectively, some methods are:

Classic Richardson $\omega \le 1/\lambda_{\max}$, C = I.

Jacobi: $\omega \ge 0$, C = D

Gauss-Seidel/ SOR: $\omega \in (0,2)$, C = D+$\omega$ L

ILU: $\omega \ge 0$, C = $L'U'$, with A = $L'U'$ + R

I also tested

Tri diagonal C = $(\omega \in (0,1))$ $\omega(D +$ first super- and sub- diagonals)

For smoothing alone these have mesh dependent convergence rates around $1 - O(h^\alpha)$, $\alpha \ge 0$, typically $\alpha \in (0,2)$. For Jacobi and Gauss Seidel $\alpha = 2$, whereas ILU and SOR will have $\alpha = 1$. These smoothing methods have the property that they smooth/ dampen high frequencies quickly. However, low frequencies are reduced ineffectively, thus resulting in the slow convergence behavior of these methods.[3,4]

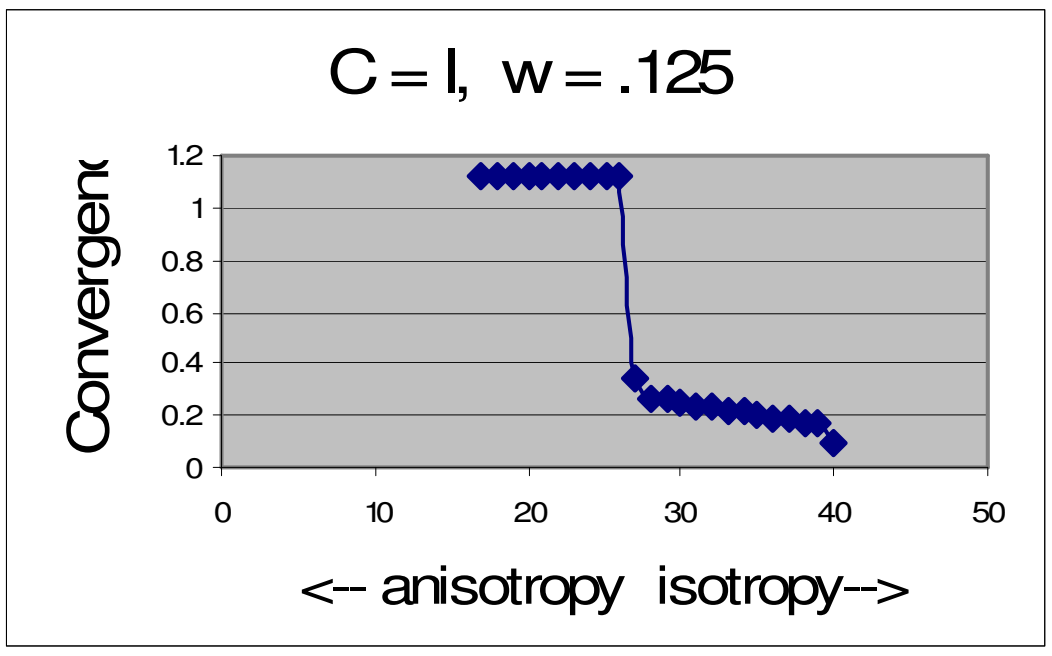

| | |
|---|---|
| Jacobi | 0 |
| TRI | 3*NEQ= number of equations |
| ADI | 2*TRI |
| GSTRI | 1*TRI |
| GSADI | 2*TRI |
| ILU | NA |

The memory required for the preconditionor used is shown in the table on the left. Before the iteration process the preconditionor memory is copied into place and sparsely inverted. Later, during smoothing the preconditionor is multiplied with the defect . For the *ADI and TRIy or GSTRIy preconditioning the defect is virtually transposed before operations and then transposed back..

**The Multigrid Technique**

The Multigrid technique uses this smoothing property[5] to dampen the low frequencies errors on coarser grids where they appear as high frequencies.

The ***Multigrid algorithm*** MG( $l$, $u_l^0$, $g_l$ )

For problem $A_l x = g$ let $l$=level, $S_l$=smoother, $P_{l-1}^l$=Prolongation, $R_l^{l-1}$=Restriction and $u_l^0$=start solution.

If $l$=1 (one level) return $A_l^{-1} g$

If $l \geq 1$:

m-pre-smoothing steps: $u_l^m = S_l u_l^{m-1} = \ldots = S_l^m u_l^0$

Course Grid Correction:

Restrict the defect: $d_l = g_l - A_l u_l^m \rightarrow g_{l-1} = R_l^{l-1} d_l$

Solve recursively: $u_{l-1}^i = MG(l-1, u_{l-1}^{i-1}, g_{l-1})$; i = 1..p, $u_{l-1}^0 = 0$

$u_l^{m+1} = u_l^m + \omega_l P_{l-1}^l u_{l-1}^p$

n-post-smoothing steps: $u_l^{m+n+1} = S_l^n u_l^{m+1}$

Finally set MG( $l$, $u_l^0$, $g_l$) $= u_l^{m+n+1}$

The $\omega_l$ parameter used above can either be a fixed value or can be chosen adaptively to minimize the error in the energy norm:

$$\omega_l = \frac{(g_l - A_l u_l^m,\ \mathrm{R}_{1\text{-}1}^{1}\mathrm{u}_{1\text{-}1})_1}{(A_l R_{1\text{-}1}^{1}\mathrm{u}_l^p,\ \mathrm{R}_{1\text{-}1}^{1}\mathrm{u}_l^p)_l}$$

$$u_{h,0} \xrightarrow[m-smoothing]{u_{h,n+1}=u_{h,n}-C^{-1}(Au_{h,n}-b)} u_{h,m} \xrightarrow[defect]{d=Au_{h,m}-b} \qquad v_h = P\tilde{u}_{2h} \xrightarrow{control} u_{h,m+1} = u_{h,m} + av_h \xrightarrow[n-postsmoothing]{} u_{h,m+n+1}$$

$$\downarrow \ restriction \qquad \uparrow \ prolongation$$

$$u_{2h} = Rd \qquad \tilde{u}_{2h}$$

$$recursion$$

$$courcest_level$$

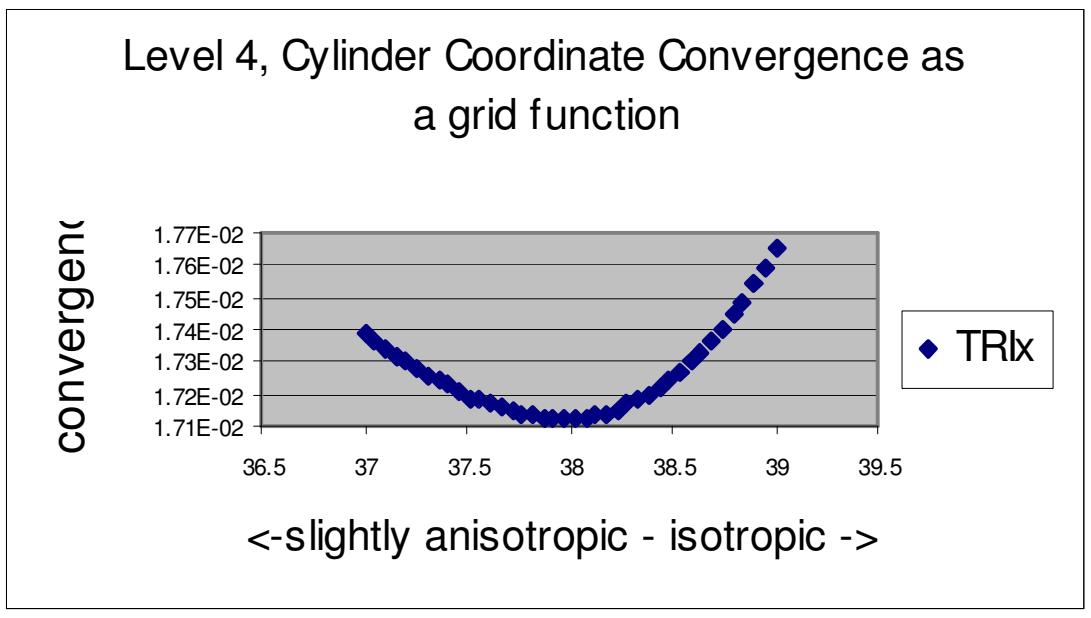

Figure: The Multigrid scheme

**Convergence studies**

**Anisotropies**

Difficulties may arise from anisotropies that are implied by operators such as $-\alpha \frac{\partial u^2}{\partial x^2} - \beta \frac{\partial u^2}{\partial y^2}$ where $\alpha \ll \beta$ or $\beta \ll \alpha$ and also high grid aspect ratios.

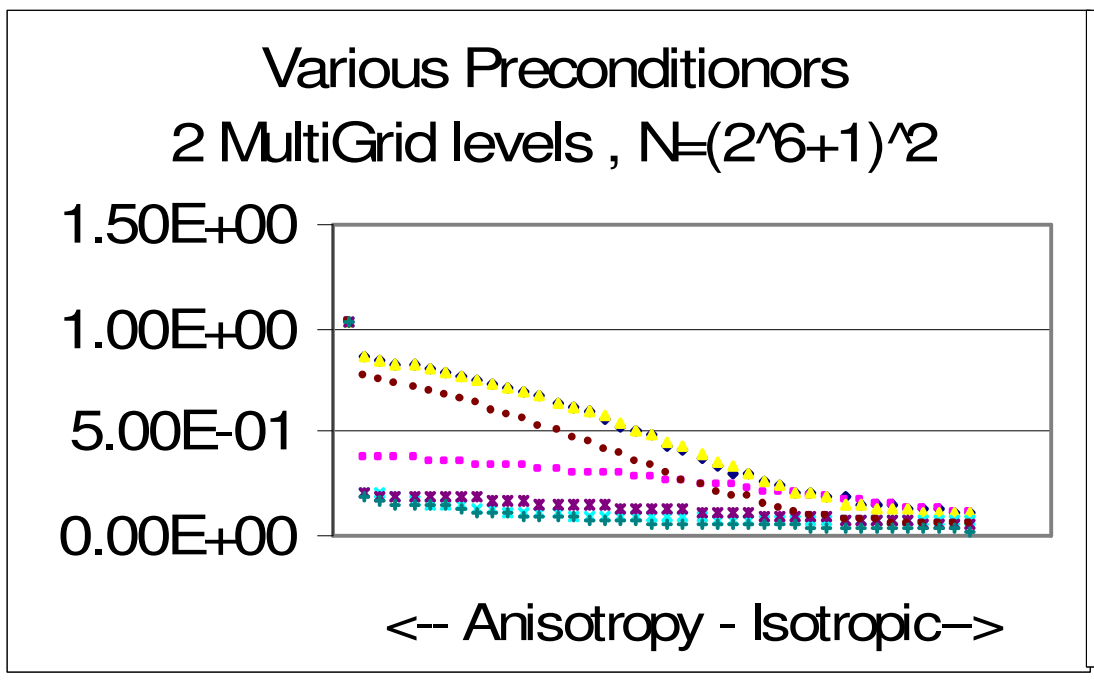

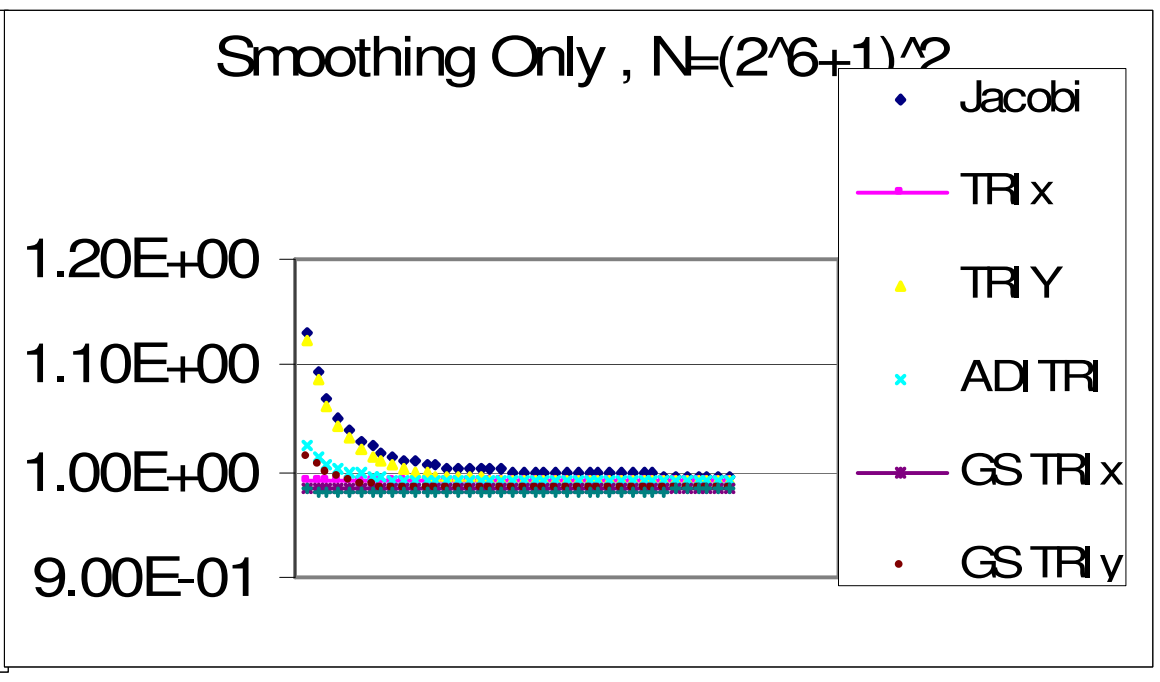

For cylinder coordinates, it is possible to find an optimal grid to improve convergence. Left Figure shows convergence rates using different preconditionors to reach a relative error of $10^{-4}$. Omega for the smoothing is .7 . Right Figure shows various smoothers on 1 level only. The slow convergence here is due to the low frequency errors and is also dependent upon the grid step size.

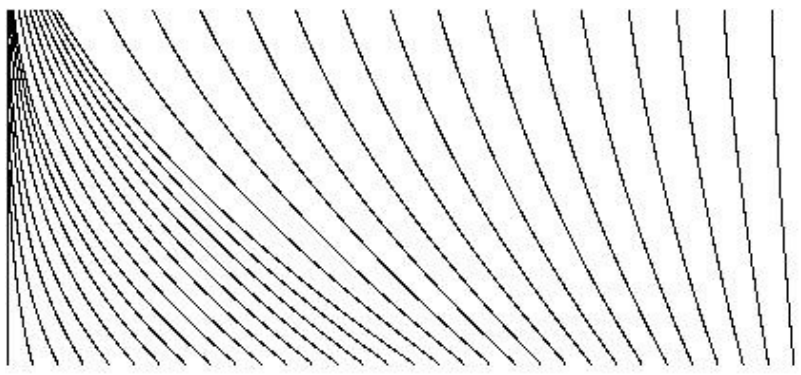

Above: Grid refinement: Bottom: equidistant, Top anisotropic (more detail in one direction). Combinations can be used to create grids with mid refinement.[6,7,8]

On the right: sample grids levels for different anisotropy parameters.

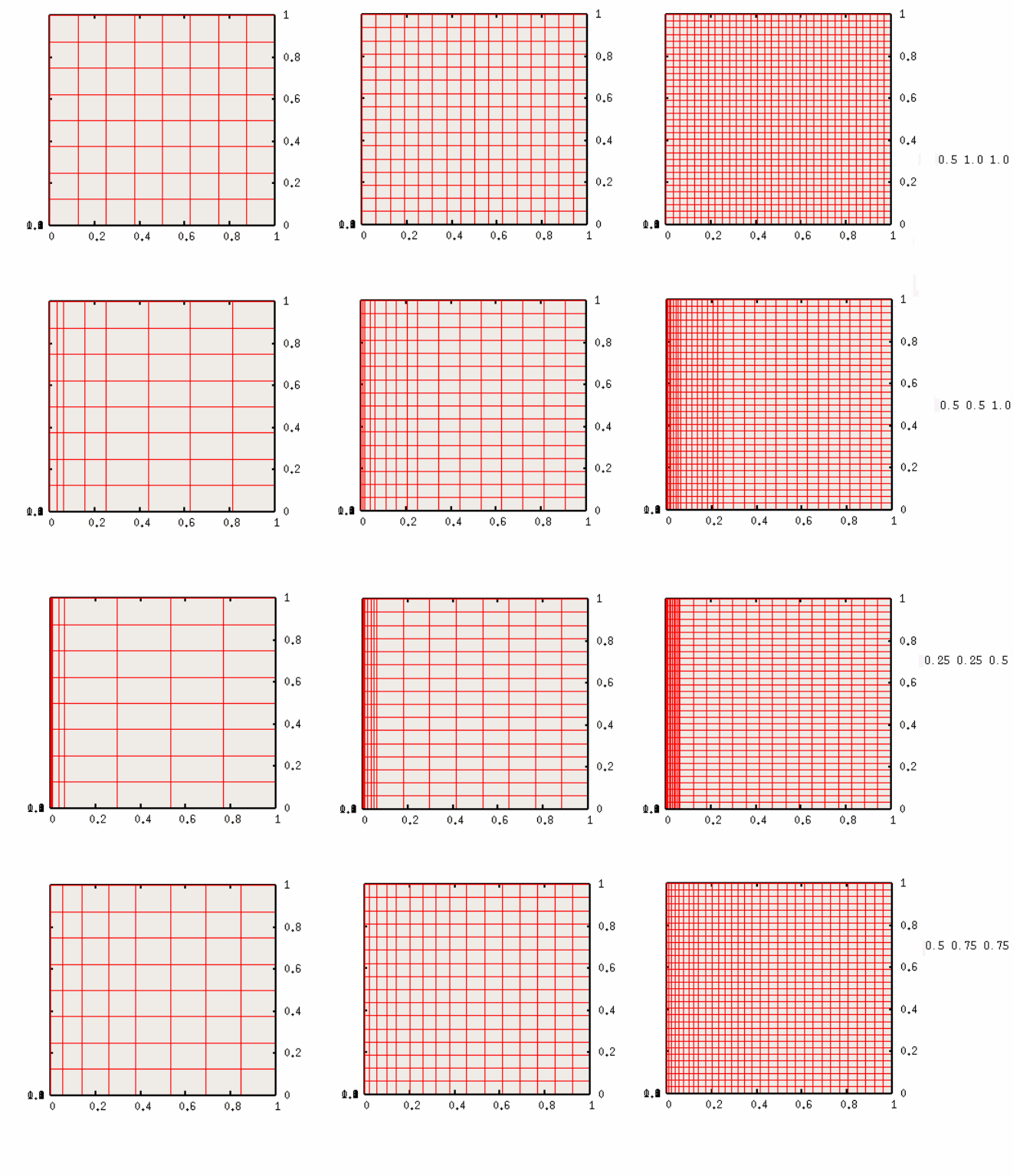

## Grid levels

One question that may arise is “Is there an optimal number of levels that could be used in order to have fast convergence while considering memory constraints? The table below shows for a sample problem on a level using a varying number of coarser grids. For this particular problem adding more than 3 grids helped none for the problem at hand.

| | |
|---|---|
| 6 Grids | 4.10541440E-03 |
| 5 Grids | 4.10529470E-03 |
| 4 Grids | 4.10518210E-03 |
| 3 Grids | 8.13862790E-03 |
| 2 Grids | 4.48421430E-01 |

**Table (left). Cartesian Isotropic Problem with Gauss Seidel Adi-preconditionor.**

## Smoothing

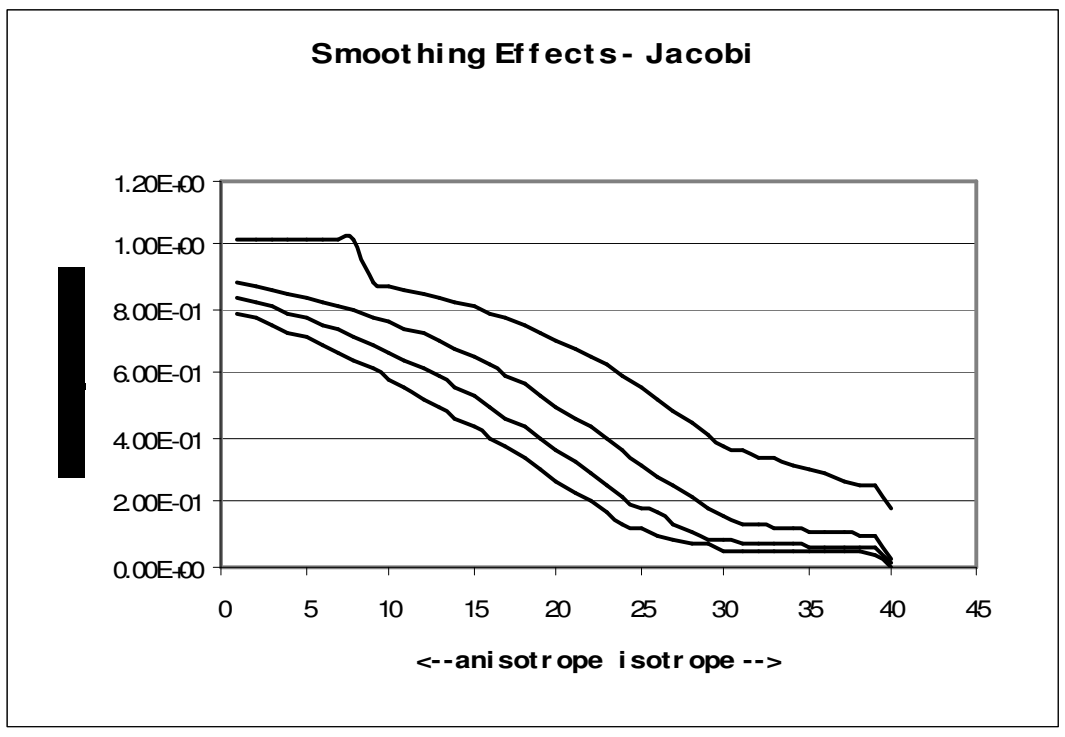

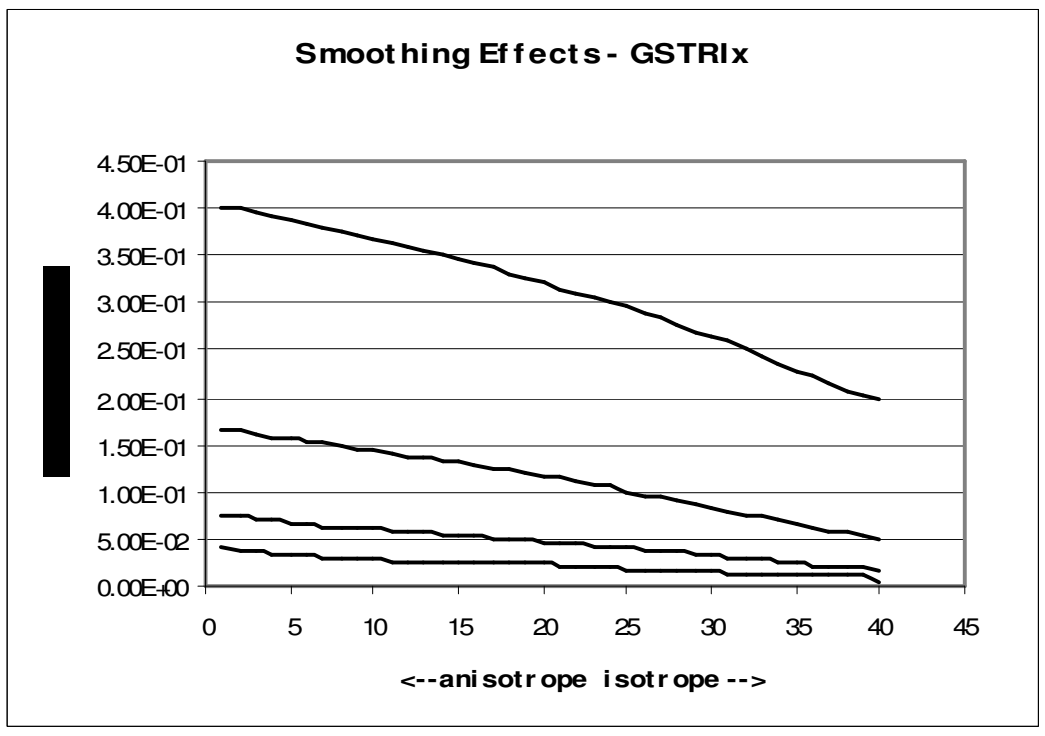

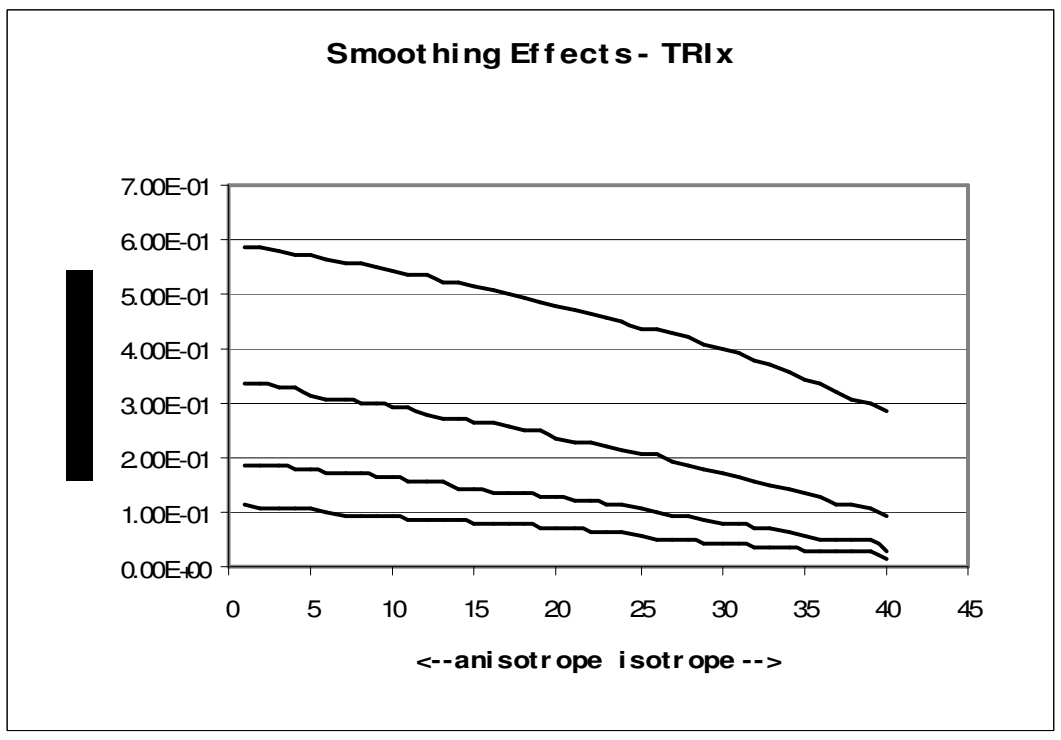

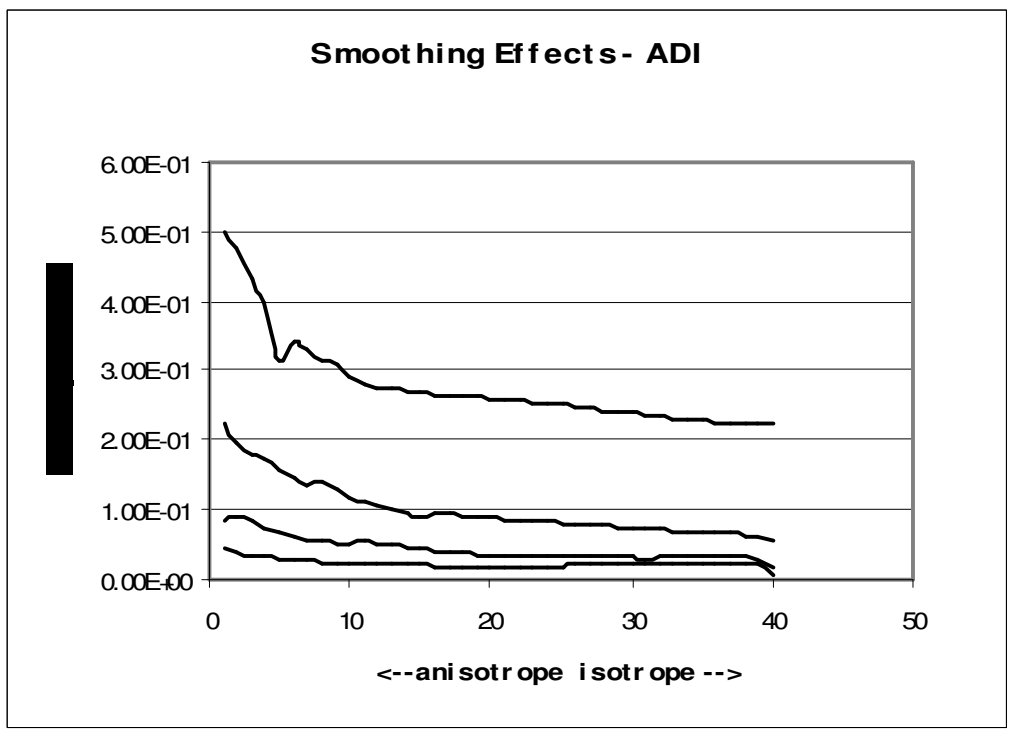

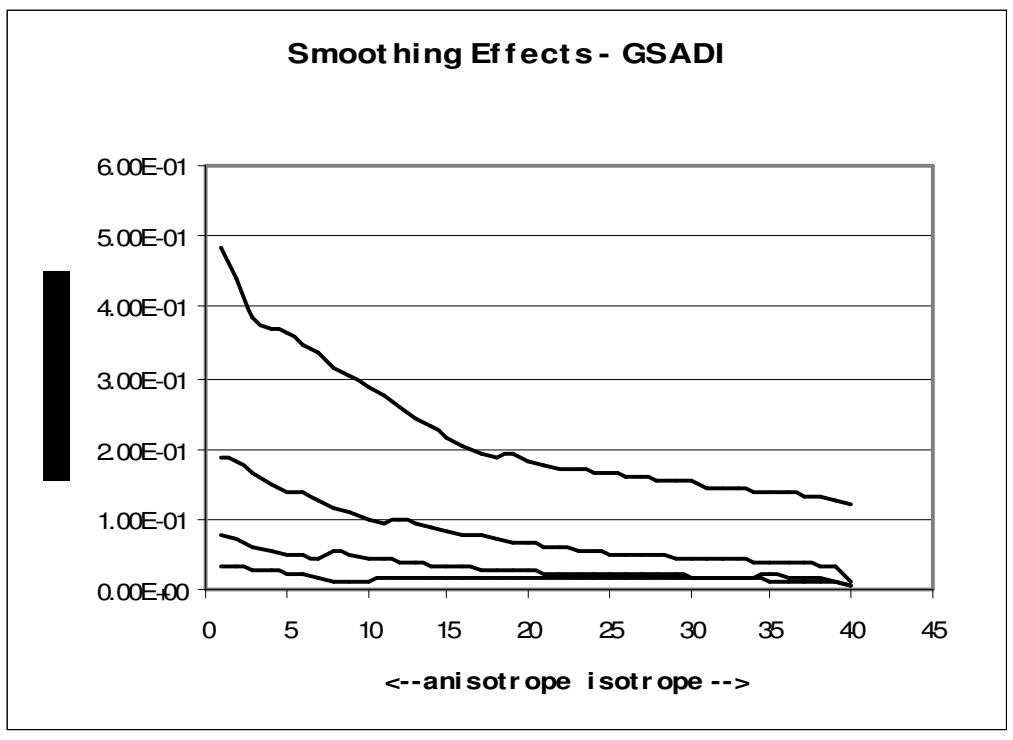

The increase in smoothing steps results in an increase of the convergence rates. However, the smoothing step is computationally wise most expensive. The figures above and to the right show 1 smoothing step (top-most curve) through 4 smoothing steps (bottom-most curve) for various preconditionors.

### Different coordinate systems

Behavior of the TRI and GSTRI preconditionors on different coordinate systems for grid refinement toward the left can be viewed on the following figures. The TRIx preconditionor faired best in Cartesian coordinates and did not converge well for spherical coordinates. The TRIy preconditionor, on the other hand did best in spherical coordinates (where coefficients of A are multiplied by $r^2\sin(\Theta)$), but did poorly in Cartesian and cylindrical coordinates with high anisotropies. The GSTRI preconditionors performed similarly but better than the TRI preconditionor (due to using more information from A. Combining the TRIx and TRIy in an ADI preconditionor brought good convergence in all coordinate systems. Finally the GSTRIx and GSTRIy preconditionors combined to GSADI resulted excellent convergence rates, but had difficulties for high anisotropies in spherical coordinates.

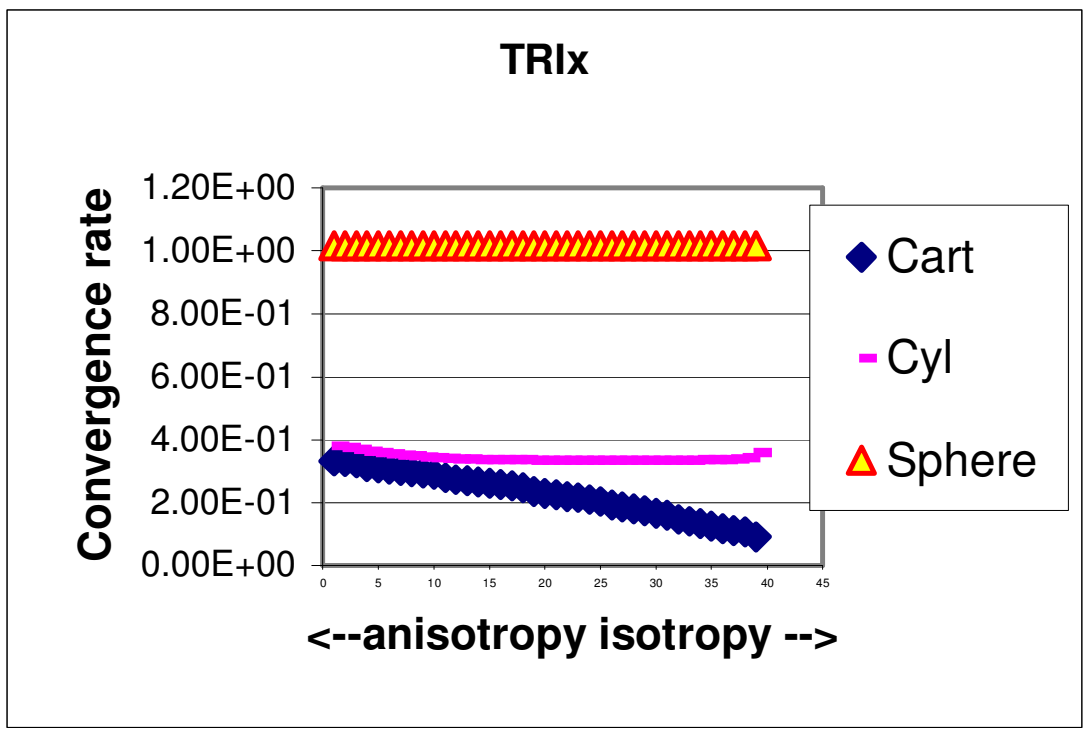

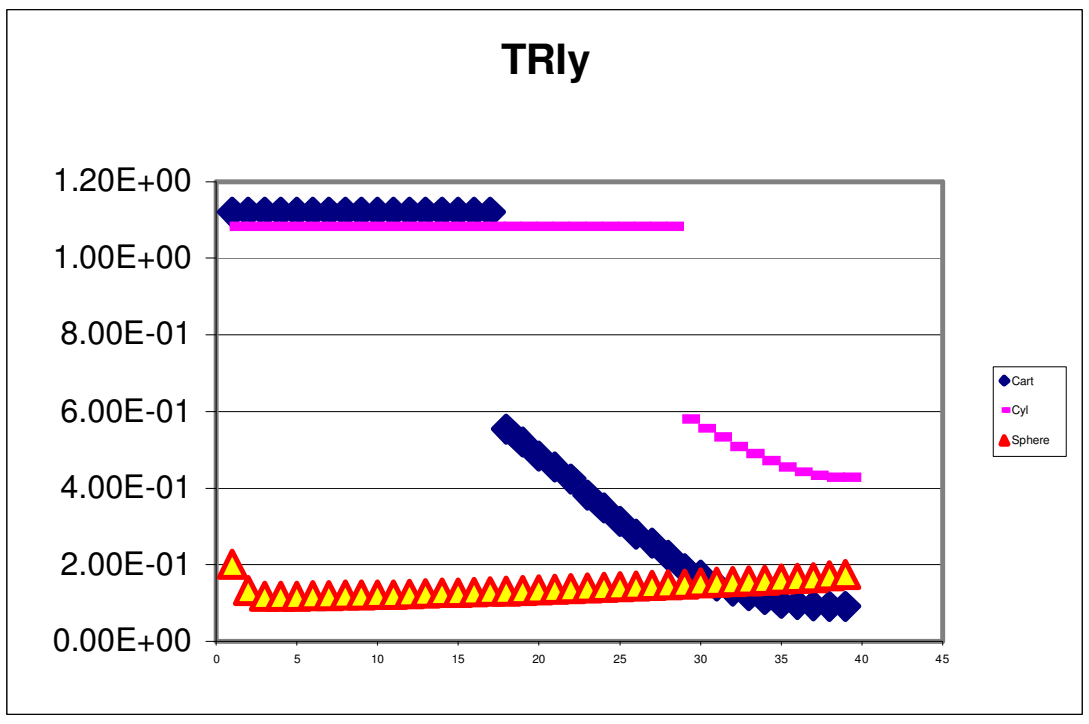

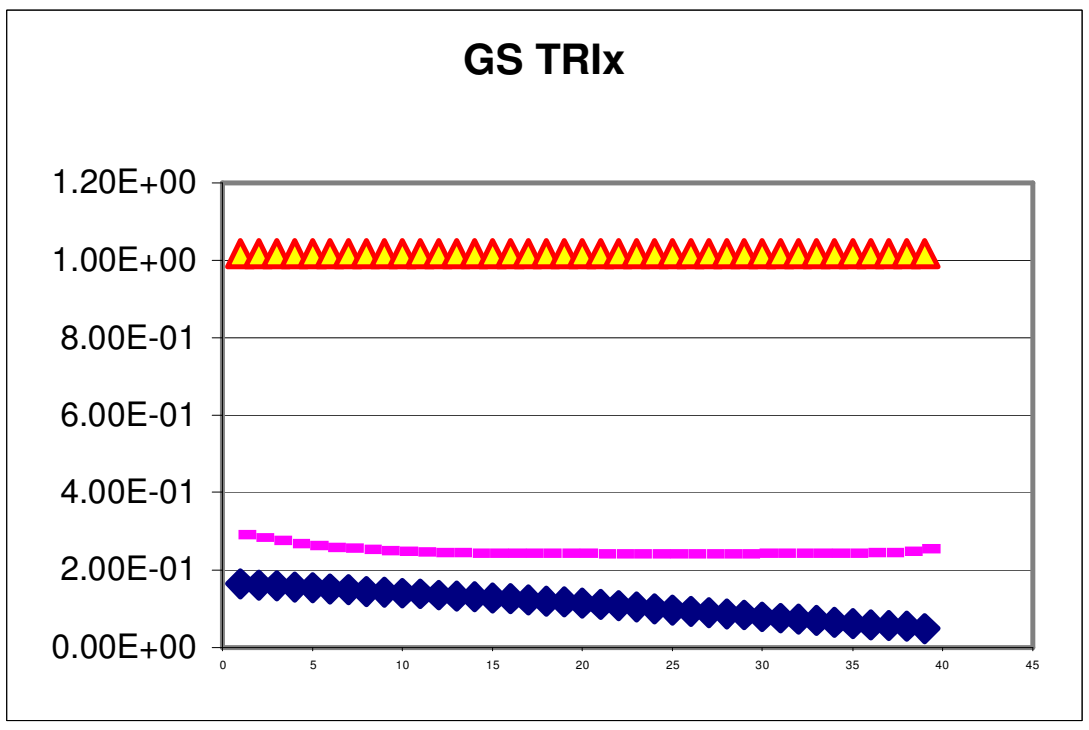

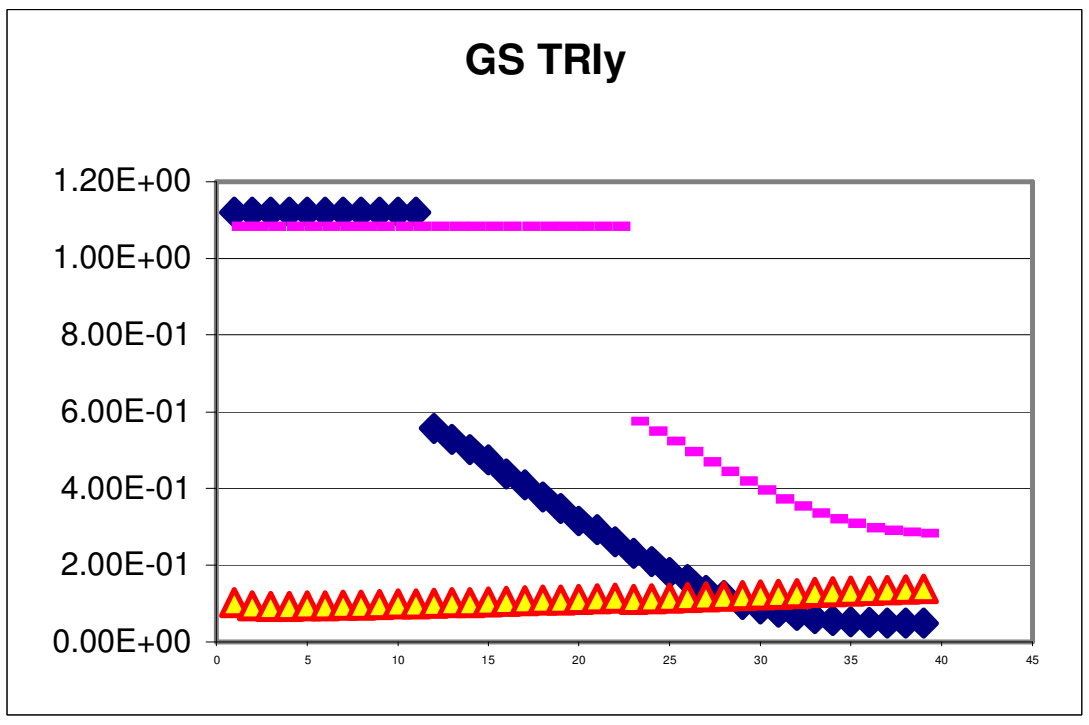

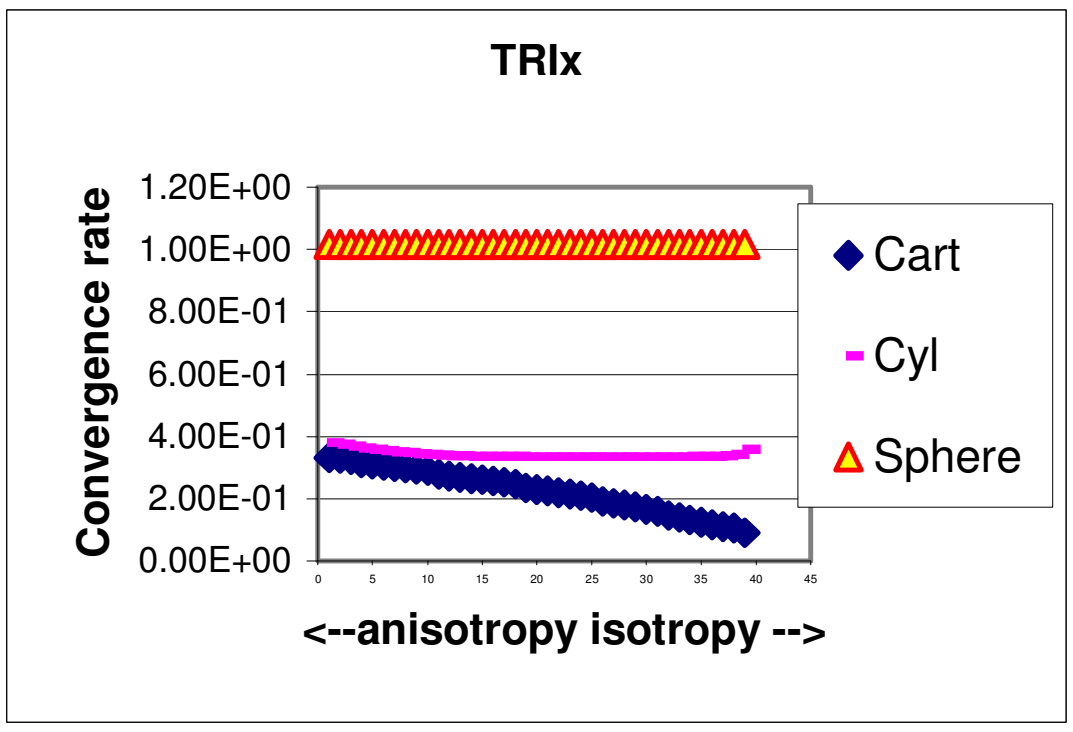

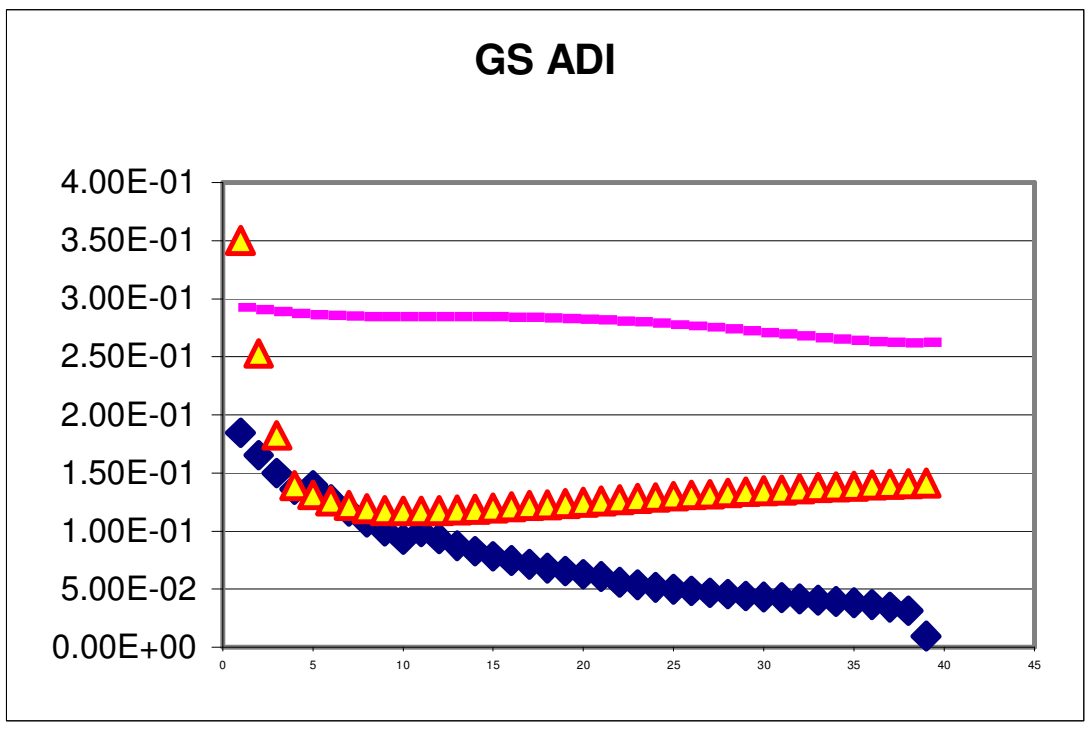

Although the TRIy preconditionor did poor for mid- to high-anisotropic grids, the combination of TRIx and TRIy preconditionors =ADI did better than just relying on a tridiagonal matrix based on direction. One be led to think that for refinement in the x-direction the tri-diagonal matrix would be best because the x-grid neighbors have more "weight" in the discretion equations. One and conclude that damping occurs only for the direction the smoother was designed for.

**Start-Vector Choice**:

The preconditionors smooth-out high frequencies well while low frequencies are smoothed on coarser grid level where they appear as high frequency errors. This leads to methods of for making better start-vectors. One way is to begin the solution method on the coarsest level and work slowly up to the finest grid level. Studies for this method can be found in [1]. Another method suggests making vector for related problems and using these as start vectors. For example, for a non-stationary problem one could find and use a solution to a stationary problem as start-vector to a non-stationary problem. Another example: In fluid dynamics, for solving problem with a very high Reynolds parameter, one can use a solution of a problem with a lower Reynolds parameter as a start vector.

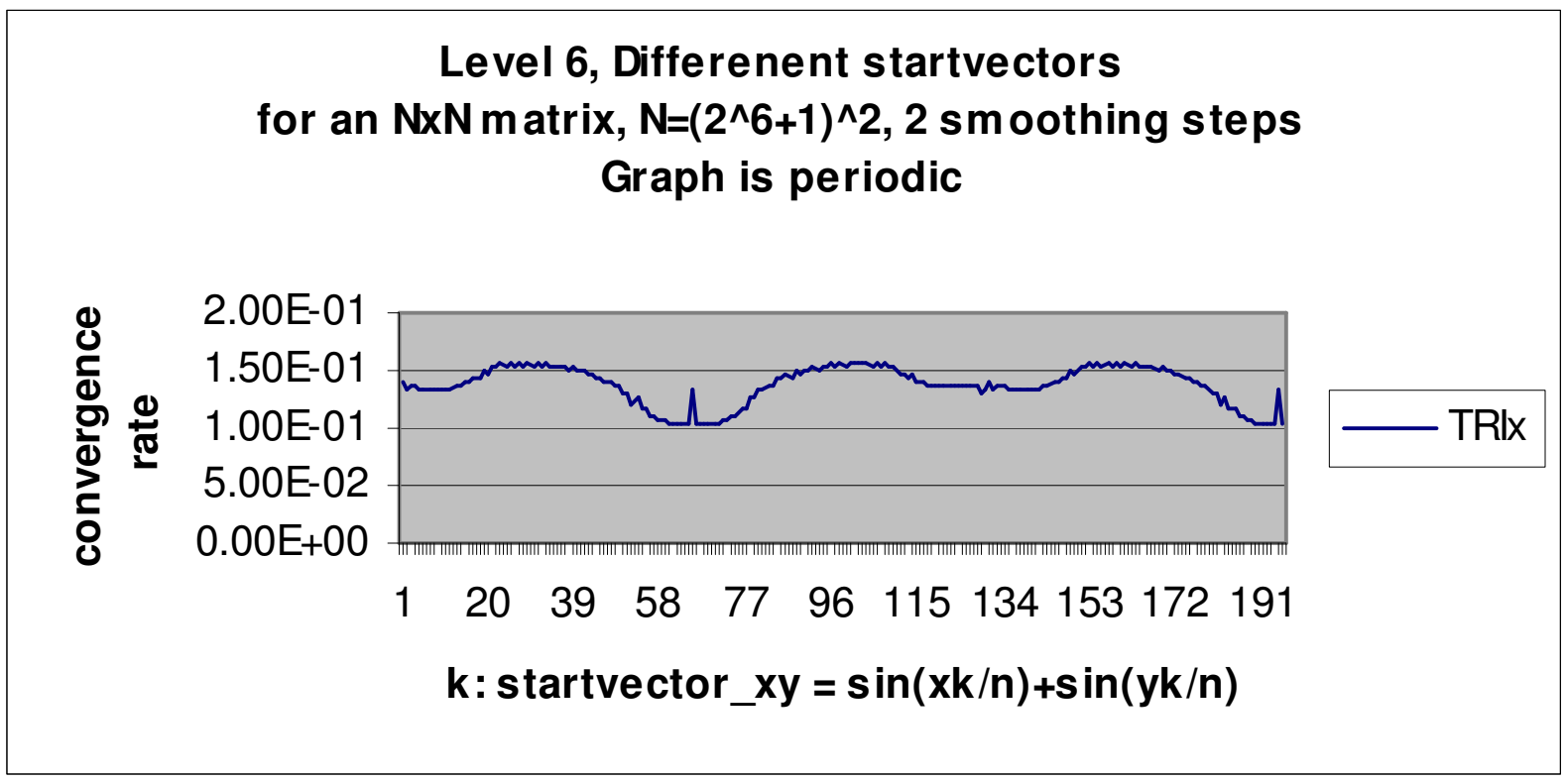

## Conclusion

In order to find an optimal problem-solving-technique while employing the Multigrid technique, it is necessary to consider the given problem, the number of unknowns, the grid and coordinate system used as well as the constraints such as the computational power and the approximation quality. Choosing a preconditionor wisely as well a good start vector choice can lead to a efficient and robust solver for a given problem.

---

[1] **Richardson, L. F.** The approximate arithmetical solution by finite differences of physical problems involving differential equations with an application to the stresses in a masonry dam. *Philos. Trans. Roy. Soc* London A 210, 307-357

[2] **L. A. Hageman and D. M. Young**, Applied iterative methods**,** Academic Press, New York, 1981: page 43

[3] **Altieri, M.; Becker, Ch.; Turek, S.**: *On the realistic performance of Linear Algebra components in iterative solvers*, Bungartz, H.-J.; Durst, F.; Zenger, Chr., Lecture Notes in Computational Science and Engineering, 8, 3-12, High Performance Scientific and Engineering Computing: Proceedings of the International FORTWIHR Conference on HPSEC, Munich, March 16-18, 1998, Springer, ISBN 3-540-65730-4, 1999

[4] **Turek, S.**: *Ein robustes und effizientes Mehrgitterverfahren zur Lösung der instationären, inkompressiblen 2-D Navier-Stokes-Gleichungen mit diskret divergenzfreien finiten Elementen*, Preprints SFB 123, Nummer 642, Universität Heidelberg, 1991

[5] **W. Hackbusch, *Multi-Grid Methods and Applications*,** Springer-Verlag**, 1985.**

[6] **Kilian, S.; Turek, S.**: *An example for parallel ScaRC and its application to the incompressible Navier-Stokes equations*, Preprints SFB 359, Nummer 98-06, Universität Heidelberg, 1998

[7] **Turek, S.; Becker, Ch.; Kilian, S.**: *Some concepts of the software package FEAST*, Palma, J. M.; Dongarra, J. J.; Hernandes, V., Lecture Notes in Computer Science, 1573, VECPAR`98 - Third International Conference for Vector and Parallel Processing, Springer, 1998

[8] **Kilian, S.**, PhD Thesis, „ScaRC, Ein verallgemeinertes Gebietszerlegungs-/Mehrgitterkonzept auf Parallelrechnern“, Dissertation, Universität Dortmund, 2001